\documentclass{article} 
\usepackage{nips15submit_e,times}
\usepackage{hyperref}
\usepackage{url}
\usepackage{booktabs}
\usepackage{amsmath,amssymb, amsthm}

\usepackage{fullpage}

\usepackage{multirow}
\usepackage{amsfonts}
\usepackage{layout}
\usepackage{amsmath}
\usepackage{amsthm}
\usepackage{amssymb}
\usepackage{url}
\usepackage{color}
\usepackage{graphicx}
\usepackage{algorithmic}
\usepackage{algorithm}
\usepackage{verbatim}
\usepackage{epsfig}
\usepackage{epstopdf} 
\usepackage{xcolor}

\PassOptionsToPackage{normalem}{ulem}
\usepackage{ulem}

\usepackage{cases}

\providecolor{added}{rgb}{0,0,1}
\providecolor{deleted}{rgb}{1,0,0}

\definecolor{orange}{RGB}{255,127,0}

\newcommand{\ve}[2]{\langle #1 ,  #2 \rangle}

\newcommand{\R}{\mathbb{R}}

\newcommand{\vc}[2]{#1^{[#2]}}

\newtheorem{assumption}{Assumption}

\theoremstyle{plain}

\theoremstyle{definition}

\newcommand*{\starnr}{\stepcounter{equation}\tag{\theequation}}

\usepackage{comment}

\title{Partitioning Data on Features or Samples in  
Communication-Efficient Distributed Optimization?}

\author{
    Chenxin Ma  \\
    Industrial and Systems Engineering\\
    Lehigh University, USA\\
    \texttt{chm514@lehigh.edu} \\
    \And
    Martin Tak\'a\v{c} \\
    Industrial and Systems Engineering\\
    Lehigh University, USA \\
    \texttt{takac.mt@gmail.com} \\
}

 \nipsfinalcopy 

\begin{document}

\maketitle

\begin{abstract}
In this paper we study the effect of the way that the data is partitioned in distributed optimization. 
The original
DiSCO algorithm 
[Communication-Efficient Distributed Optimization of
Self-Concordant Empirical Loss, Yuchen Zhang and
Lin Xiao, 2015]
partitions the input data based on samples. 
We describe how the original algorithm has to be modified to allow partitioning on features and show its efficiency both in theory and also in practice.
\end{abstract}

\section{Introduction}

As the size of the datasets becomes larger and larger, distributed optimization methods for machine learning have become increasingly important \cite{bertsekas1989parallel,dekel2012optimal,shamir2014distributed}. Existing mehods often require a large amount of communication between computing nodes \cite{yang2013trading,jaggi2014communication,ma2015adding,yang2013analysis}, which is typically several magnitudes slower than reading data from their own memory \cite{marecek2014distributed}. Thus, distributed machine learning suffers from the communication bottleneck on real world applications.

In this paper we focus on the regularized empirical risk minimization problem. Suppose we have $n$   data samples $\{x_i, y_i\}_{i=1}^n$,
where each $x_i \in \R^d$ (i.e. we have $d$ features), $y_i \in \R$.
We will denote by  $X := [x_1,...,x_n]\in\R^{d\times n}$. The optimization problem is to minimize  the regularized empirical loss (ERM)
\begin{equation}\label{eq:pro1}
  f(w) := \frac{1}{n} \sum_{i=1}^n \phi_i(w,x_i) + \frac{\lambda}{2} \|w\|_2^2,
\end{equation}
where the 
first part is the \emph{data fitting term}, 
 $\phi: \R^d\times \R^d \rightarrow \R$ is a loss function which typically depends on $y_i$.
Some popular loss functions includes hinge loss $\phi_i(w, x_i)=\max\{0, 1-y_i w^T x_i\}$, square loss
  $\phi_i(w, x_i)=(y_i- w^T x_i)^2$
  or logistic loss
  $\phi_i(w, x_i)=\log(1+\exp( -y_i w^T x_i))$. 
 The second part of objective function \eqref{eq:pro1} is $\ell_2$ regularizer ($\lambda>0$)
 which helps to prevent over-fitting of the data.
 
We assume that the loss function $\phi_i$ is convex and  self-concordant \cite{zhang2015communication}:
\begin{assumption}\label{ass:selfc}
For all $ i \in [n]:=\{1,2,\dots,n\}$ the convex function $\phi$ is self-concordant with parameter $M$ i.e. the following inequality holds:
\begin{equation}
  |u^T (f'''(w)[u])u| \leq M(u^T f''(w) u) ^{\frac{3}{2}}
\end{equation}
for any $u\in\R^d$ and $w\in dom(f)$, where $f'''(w)[u]:= \lim_{t\rightarrow 0} \frac{1}{t} \large(f''(w+tu)- f''(w)\large)$.
\end{assumption}

There has been an enormous interest in large-scale machine learning problems and many parallel \cite{bradley2011parallel,recht2011hogwild}
or distributed algorithms have been proposed 
\cite{agarwal2011distributed,takavc2015distributed,
richtarik2013distributed,shamir2013communicationNewton,
lee2015distributed}. 

The main bottleneck in distributed computing --communication-- was
handled by many researches differently. Some work considered ADMM type methods \cite{boyd2011distributed,deng2012global},
another used block-coordinate type algorithms
\cite{lee2015distributed,yang2013trading,jaggi2014communication,ma2015adding}, where they tried to solve the local sub-problems more accurately (which should decrease the overall communications requirements  when compared with more basic approaches \cite{takavc2013mini,takavc2015distributed}).

\section{Algorithm}
\begin{algorithm}[t]
\caption{High-level DiSCO algorithm}
\label{Disco}
\begin{algorithmic}[1]
\STATE\textbf{Input: } parameters $\rho, \mu \geq 0$, number of iterations $K$
\STATE Initializing $w_0$.
\FOR  {k = 0,1,2,...,K}
\STATE Option 1: Given $w_k$, run DiSCO-S PCG Algorithm \ref{A1}, get $v_k$ and $\delta_k$
\STATE Option 2: Given $w_k$, run DiSCO-F PCG Algorithm \ref{A2}, get $v_k$ and $\delta_k$
\STATE Update $w_{k+1} = w_k - \frac{1}{1+ \delta_k} v_k$
\ENDFOR
\STATE \textbf{Output: $w_{K+1}$}
\end{algorithmic}
\end{algorithm}
We assume that we have  $m$ machines (computing nodes) available
which can communicate between each other over the network.
We assume that the space needed to store the data matrix $X$ exceeds the memory of every single node. Thus we have to split the data (matrix $X$) over the $m$ nodes. The natural question is: \emph{How to split the data into $m$ parts?}
There are many possible ways, but two obvious ones:
\begin{enumerate}
\item split the data matrix $X$ by rows (i.e.  create $m$ blocks by rows);
Because rows of $X$ corresponds to features, we will denote the algorithm which is using this type of partitioning as \emph{DiSCO-F};

\item split the data matrix $X$ by columns;
Let us note that columns of $X$ corresponds to samples  we will denote the algorithm which is using this type of partitioning as \emph{DiSCO-S};
\end{enumerate}
 Notice that the DiSCO-S is exactly the same as DiSCO proposed and analyzed 
 in \cite{zhang2015communication}. 
 In each iteration of Algorithm \ref{Disco}, wee need to compute an inexact Newton step $v_k$ such that $\|f''(w_k)v_k - f'(w_k)\|_2 \leq \epsilon_k$, which is an approximate solution to the Newton system $f''(w_k)v_k = \nabla f (w_k)$. The discussion about how to choose  $\epsilon_k$ and $K$ and a convergence guarantees  for  Algorithm \ref{Disco} can be found   in \cite{zhang2015communication}.
The main goal of this work is to analyze the algorithmic modifications to DiSCO-S when the partitioning type is changed. It will turn out that partitioning on features (DiSCO-F) can lead to algorithm which uses less communications (depending on the relations between $d$ and $n$) (see Section \ref{sec:numExperiments}).

 \begin{algorithm}[b]
\caption{Distributed DiSCO-S: PCG algorithm -- data  partitioned by  samples}
\label{A1}
\begin{algorithmic}[1]
\STATE\textbf{Input: } $w_k\in \R^{d}$, and $\mu \geq 0$.\hfill communication (Broadcast $w_k\in\R^d$ and reduceAll $\nabla f_i(w_k)\in\R^d$ )
\STATE\textbf{Initialization: } Let $P$ be computed as \eqref{eq:precon}. $v_0 = 0$, $s_0 = P^{-1} r_0$, $r_0 =\nabla f(w_k)$, $u_0 = s_0$.
\FOR {$t= 0,1 ,2,...$}
\STATE Compute $Hu_t$       \hfill communication (Broadcast $u_t\in\R^d$ and reduceAll $f{''}_i(w_k)u_t\in\R^d$ )
\STATE {Compute $\alpha_t = \frac{ \ve{r_t}{s_t} }{\ve{u_t}{Hu_t} }$}  
\STATE Update $v_{(t+1)} = v_t+ \alpha_t u_t, Hv_{(t+1)} = Hv_t + \alpha_t Hu_t, r_{t+1} = r_t - \alpha_t H u_t$.
\STATE Update $s_{(t+1)} = P^{-1} r_{(t+1)}$.
\STATE {Compute $\beta_t = \frac{  \ve{r_{(t+1)} }{ s_{(t+1)} }}{ \ve{r_t}{s_t} }$}
\STATE Update  $u_{(t+1)} = s_{(t+1)} + \beta_t u_t$.
\STATE\textbf{until: $\|r_{(r+1)}\|_2 \leq \epsilon_k$}
\ENDFOR
\STATE \textbf{Return: $v_k =v_{t+1}$, $\delta_k = \sqrt{v_{(t+1)}^THv_{t} + \alpha_t v_{(t+1)}^T Hu_t}$}
\end{algorithmic}
\end{algorithm}
\paragraph{DiSCO-S Algorithm.} 
If the dataset is partitioned by samples, such that $j$--th node will only store $X_j =[x_{j,1},...,x_{j,n_j}] \in\R^{d\times n_j}$, which is a part of $X$, then each machine can evaluate a local empirical loss function
\begin{equation}
   f_j(w):= \frac{1}{n_j} \sum_{i=1}^{n_j} \phi(w,x_{j,i}) + \frac{\lambda}{2} \|w\|_2^2.
\end{equation} 
Because $\{X_j\}$ is a partition of $X$ we have 
  $\sum_{j=1}^m n_j = n$, our goal now becomes to minimize the function $f(w) = \frac{1}{m} \sum_{h=1}^{m} f_j(w)$. Let $H$ denote the Hessian $f''(w_k)$.
For simplicity in this paper we   consider only square loss  and hence in this case $f''(w_k)$ is constant (independent on $w_k$).

In Algorithm \ref{A1}, each machine will use its local data to compute the local gradient and local Hessian and then aggregate them together. We also have to choose one machine as the master, which computes all the vector operations of PCG loops (Preconditioned Conjugate Gradient), i.e., step 5-9 in Algorithm \ref{A1}. 

The preconditioning matrix for PCG is defined only on master node and consists of the local Hessian approximated by a subset of data available on master node with size $\tau$, i.e.
\begin{equation}\label{eq:precon}
  P =   \frac{1}{\tau} \sum_{j=1}^\tau \phi^{''} (w, x_{1,j}) + \mu I,
\end{equation}
where $\mu$ is a small regularization parameter. Algorithm \ref{A1} presents the distributed PCG mathod for solving the preconditioning linear system
\begin{equation}
  P^{-1} H v_k = P^{-1} \nabla f(w_k).
\end{equation}

\paragraph{DiSCO-F Algorithm.} 
If the dataset is partitioned by features, then $j$th machine will store $X_j = [\vc{a_1}{j},...,\vc{a_n}{j}]\in\R^{d_j\times n}$, which contains all the samples, but only with a subset of features. Also, each machine will only store $\vc{w_k}{j}\in\R^{d_j}$ and thus only be responsible for the computation and updates of $\R^{d_j}$ vectors. By doing so, we only need one ReduceAll on a vector of length $n$, in addition to two ReduceAll on scalars number.

\begin{algorithm}[h]
\caption{Distributed DiSCO-F: PCG algorithm -- data  partitioned by features}
\label{A2}
\begin{algorithmic}[1]
\STATE\textbf{Input: } $\vc{w_k}{i}\in \R^{d_i}$ for $i = 1,2,...,m$, and $\mu \geq 0$. 
\STATE\textbf{Initialization: } Let $P$ be computed as \eqref{eq:precon}. $\vc{v_0}{i} = 0$, $\vc{s_0}{i} = {(P^{-1})}^{[i]} \vc{r_0}{i}$, $\vc{r_0}{i} = f'(\vc{w_k}{i})$, $\vc{u_0}{i} = \vc{s_0}{i}$.
\WHILE {$\|r_{r+1}\|_2 \leq \epsilon_k$}
\STATE Compute $\vc{(Hu_t)}{i}$. \hfill communication (ReduceAll an $\R^{d_i}$ vector)
\STATE {Compute $\alpha_t = \frac{\sum_{i=1}^{m}  \ve{\vc{r_t}{i}}{\vc{s_t}{i}} }{\sum_{i=1}^{m}  \ve{\vc{u_t}{i}}{\vc{(Hu_t)}{i}} }$}. \hfill communication (ReduceAll a number)
\STATE Update $\vc{v_{t+1}}{i} = \vc{v_t}{i} + \alpha_t \vc{u_t}{i}, \vc{(Hv_{t+1})}{i} = \vc{(Hv_t)}{i} + \alpha_t \vc{(Hu_t)}{i}, \vc{r_{t+1}}{i} = \vc{r_t}{i} - \alpha_t \vc{(Hu_t)}{i}$.
\STATE Update $\vc{s_{t+1}}{i} = (P^{-1})^{[i]} \vc{r_{t+1}}{i}$.
\STATE {Compute $\beta_t = \frac{\sum_{i=1}^{m}  \ve{\vc{r_{t+1}}{i}}{\vc{s_{t+1}}{i}} }{\sum_{i=1}^{m}  \ve{\vc{r_t}{i}}{\vc{s_t}{i}} }$}. \hfill communication (ReduceAll a number)
\STATE Update  $\vc{u_{t+1}}{i} = \vc{s_{t+1}}{i} + \beta_t \vc{u_t}{i}$.
\ENDWHILE
\STATE Compute $\vc{\delta_k}{i}= \sqrt{{\vc{v_{t+1}}{i}}^T(Hv_t)^{[i]} + \alpha_t{\vc{v_{t+1}}{i}}^T (Hu_t)^{[i]}}$.
\STATE \textbf{Integration:  $v_{k}=[\vc{v_{t+1}}{1},...,\vc{v_{t+1}}{m}]$, $\delta_{k}=[\vc{\delta_{t+1}}{1},...,\vc{\delta_{t+1}}{m}]$} \hfill communication (Reduce an $\R^{d_i}$ vector)
\STATE \textbf{Return: $v_{k}$, $\delta_{k}$}
\end{algorithmic}
\end{algorithm}

\paragraph{Comparison of Communication and Computational Cost.}  
In Table \ref{tab:addlabel} we compare the communication cost for the two approaches DiSCO-S/DiSCO-F. As it is obvious from the table, DiSCO-F 
requires only one reduceAll of a vector of length $n$, whereas the DiSCO-S needs one reduceAll of a vector of length $d$ and one broadcast of vector of size $d$. So roughly speaking, when $n < d$ then DiSCO-F will need less communication.
However, very interestingly, the advantage of DiSCO-F is the fact that it uses CPU on every node more effectively. It also requires less total amount of work to be performed on each node, leading to more balanced and efficient utilization of nodes.

\begin{table}
 \caption{Comparison of computation and communication between different ways of partition on data.}
  \label{tab:addlabel}%
  \centering
    \begin{tabular}{|c|c|r|c|c|}
    \hline
    \multicolumn{3}{|c|}{} & partition by samples & partition by features \\
    \hline
    \multirow{8}[0]{*}{computation} &\multirow{4}[0]{*}{master} & matrix-vector multiplication & $ 1 (\R^{d\times d}\times \R^d)$ &$1 (\R^{d_1\times d_1}\times \R^{d_1})$\\
          &   &back solving linear system & 1 ($\R^d$)  & 1 $(\R^{d_1})$ \\
          &   &sum of vectors & 4 ($\R^d$)  & 4 ($\R^{d_1}$)\\
          &   &inner product of vectors & 4 ($\R^{d}$) & 4 ($\R^{d_1}$)\\
         \cline{2-5} 
         & {\multirow{4}[0]{*}{nodes}} & matrix-vector multiplication & 1 $ (\R^{d\times d}\times \R^d)$& $1 (\R^{d_1\times d_i}\times \R^{d_i})$ \\
          &  & back solving linear system    &   0 & 1 $(\R^{d_i})$ \\
          &   & sum of vectors      &0  &  4 $(\R^{d_i})$\\
          &   &  inner product of vectors     & 0 & 4 $(\R^{d_i})$ \\
  \hline  
  \multirow{2}[0]{*}{communication} & \multicolumn{2}{|c|}{Broadcast} &  one $\R^d$ vector & 0  \\
           \cline{2-5} 
        & \multicolumn{2}{|c|}{ReduceAll}  &  one $\R^d$ vector & one $\R^n$ vector, 2 $\R^1$ \\
    \hline
    \end{tabular}%
 
\end{table}%

\section{Numerical Experiments}  
  \label{sec:numExperiments}
 
\begin{figure}
\centering

\includegraphics[scale=.18]{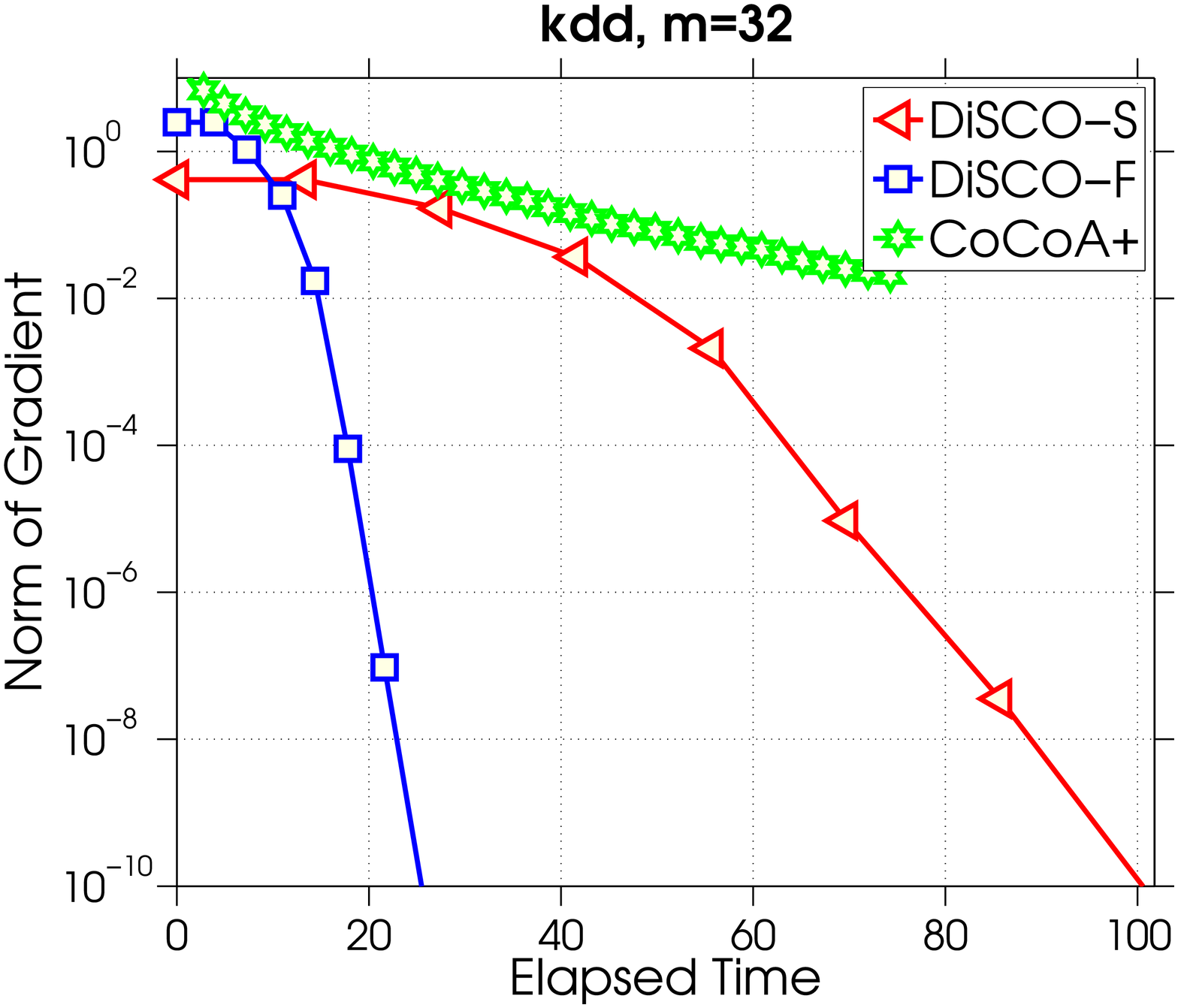}
\includegraphics[scale=.18]{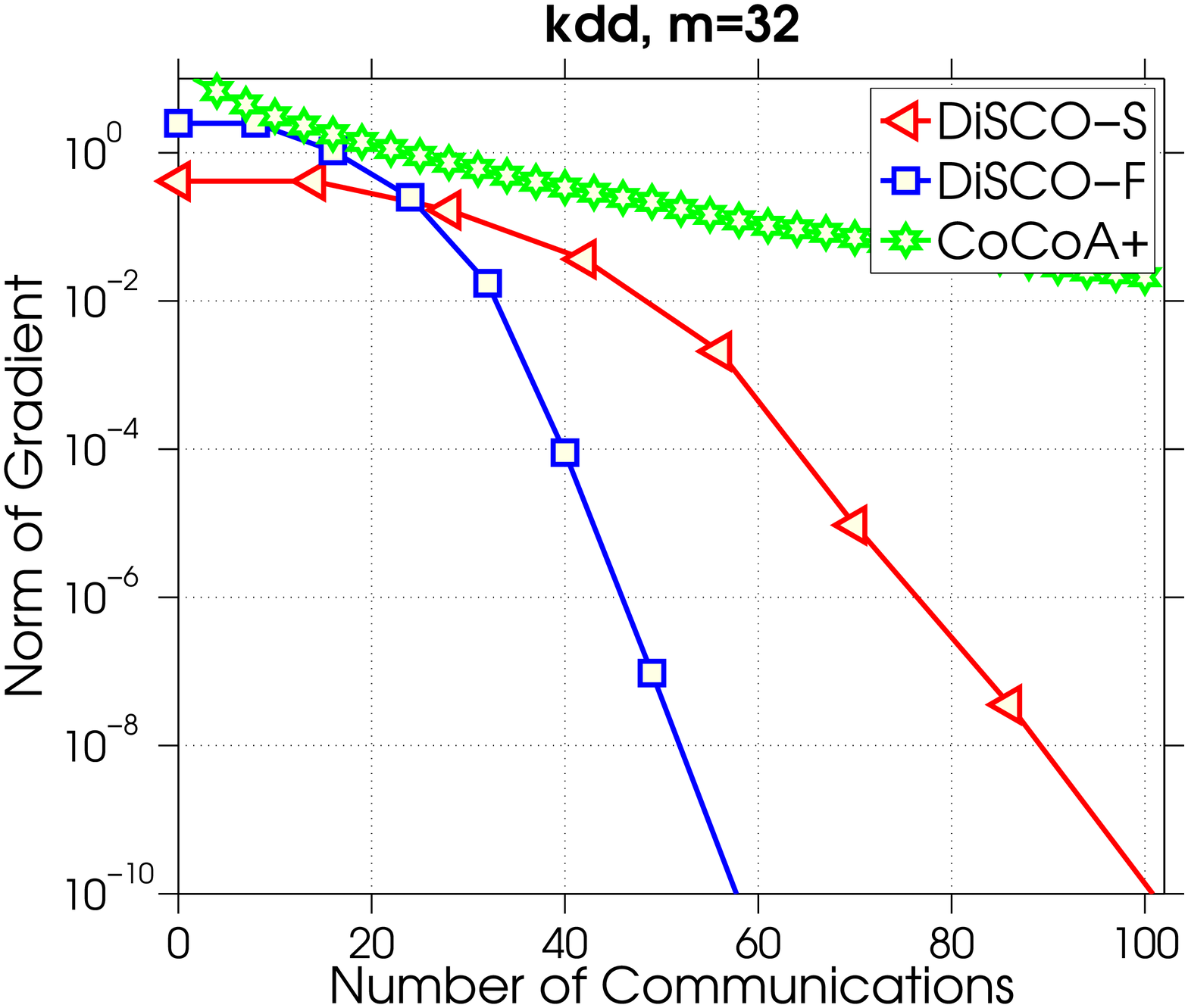}
\includegraphics[scale=.18]{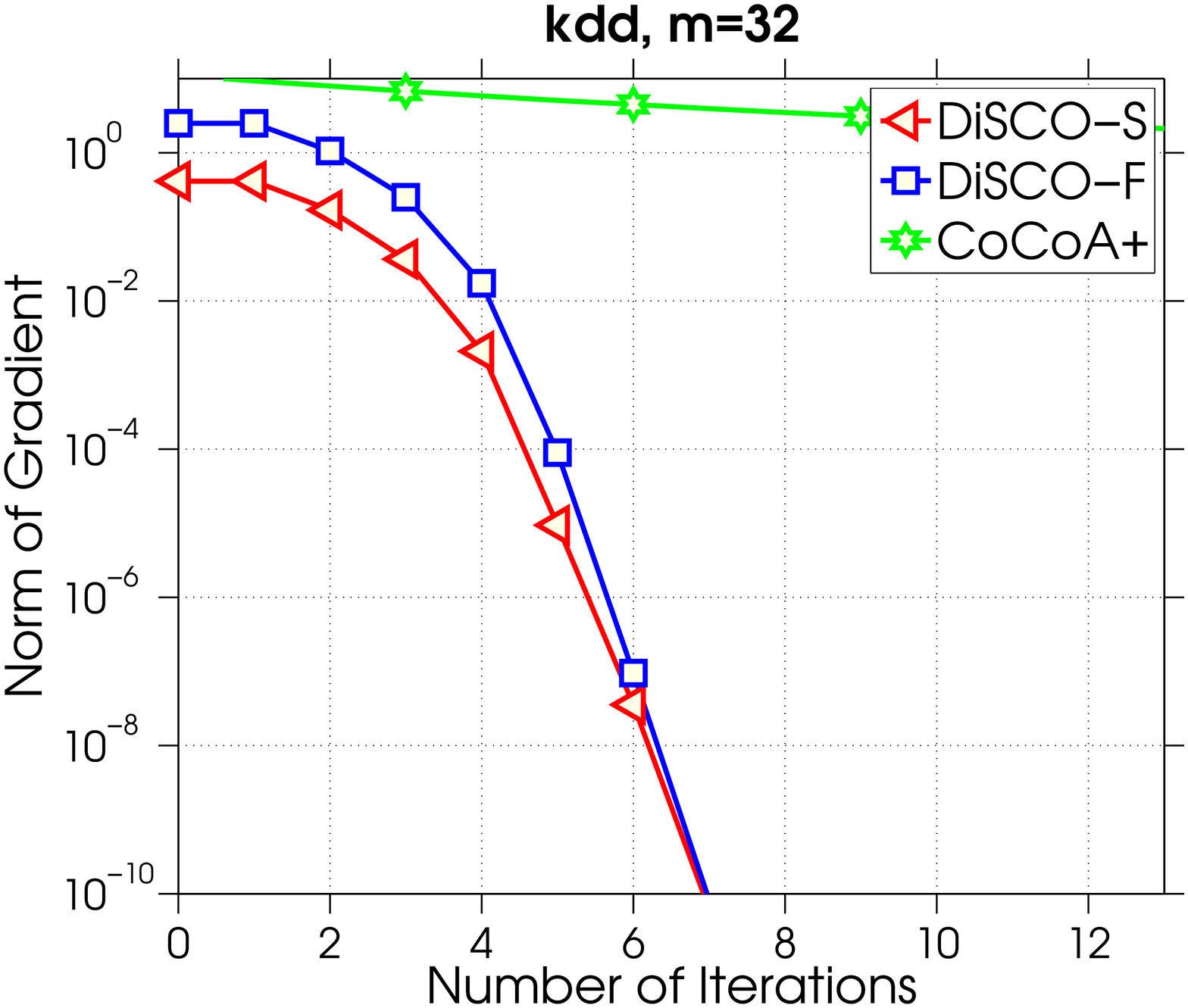}

\includegraphics[scale=.18]{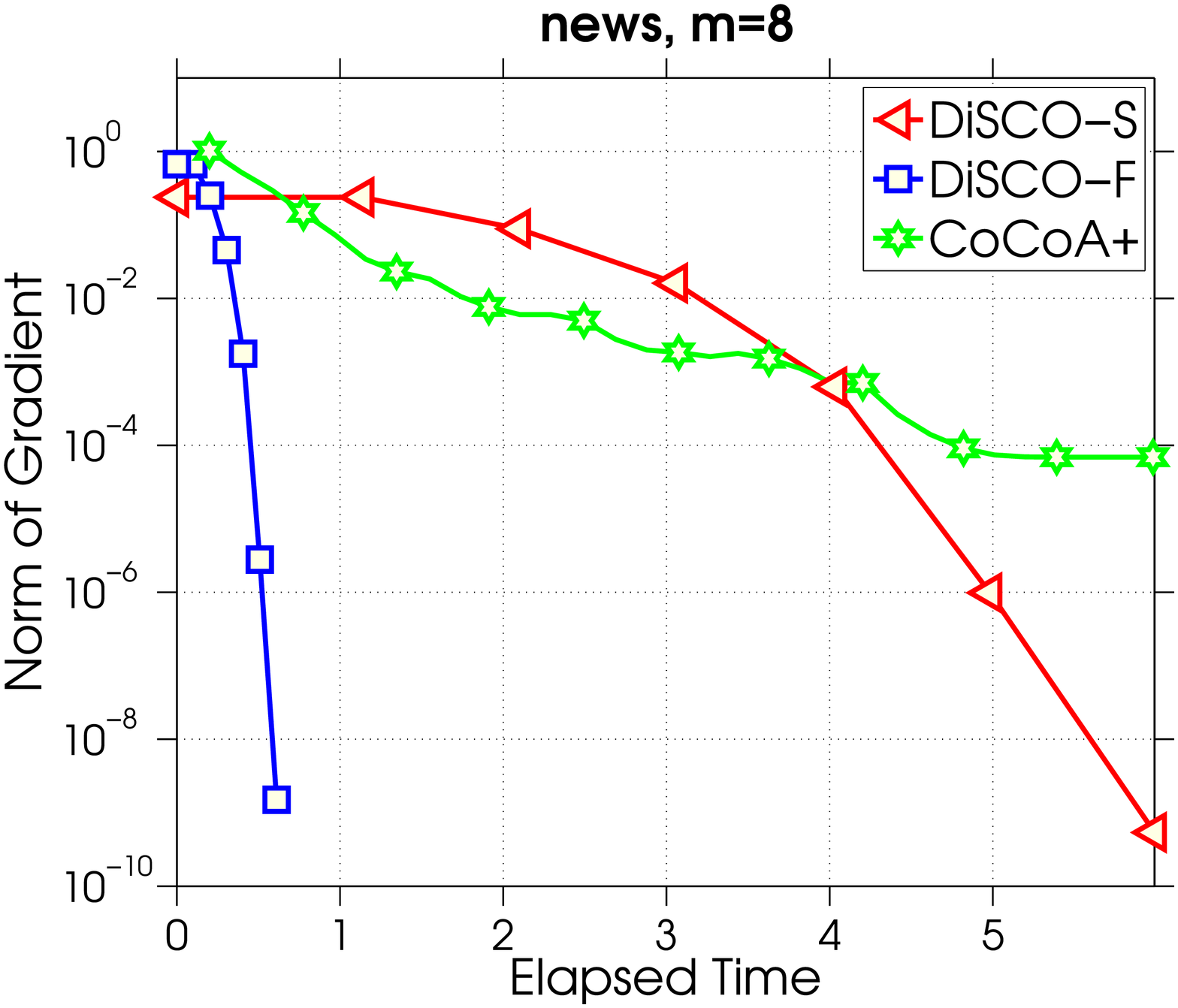}
\includegraphics[scale=.18]{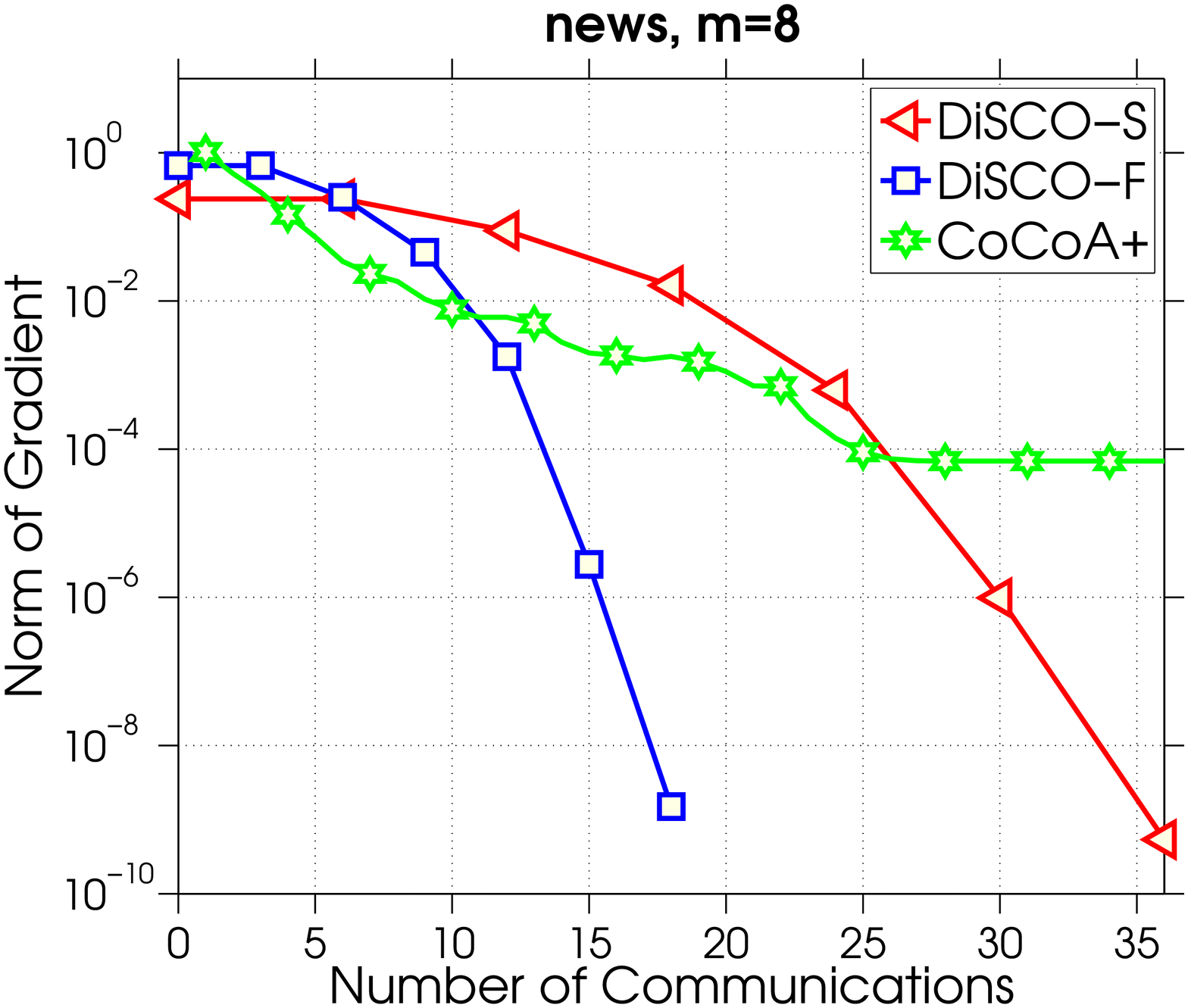}
\includegraphics[scale=.18]{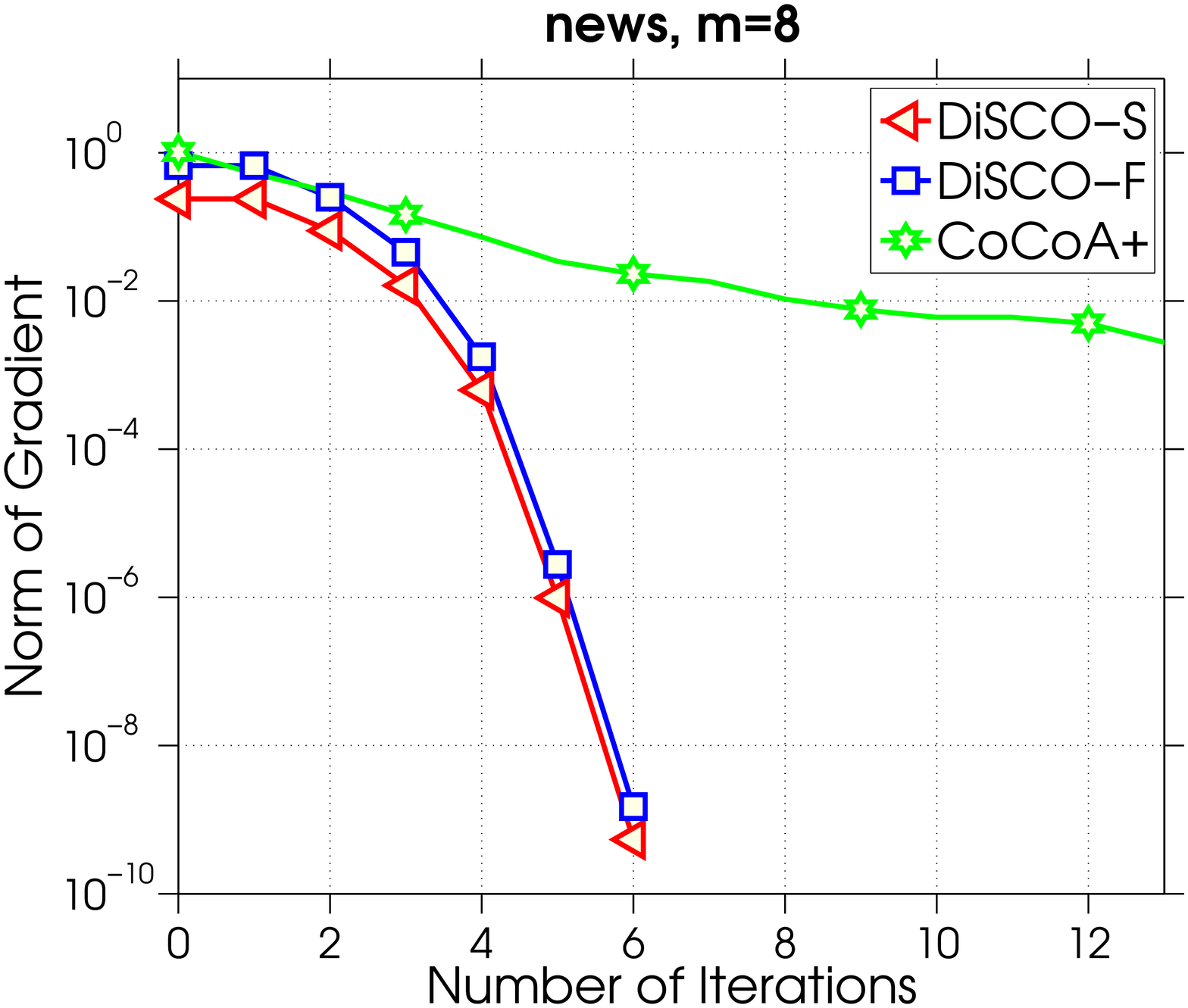}
 
\includegraphics[scale=.18]{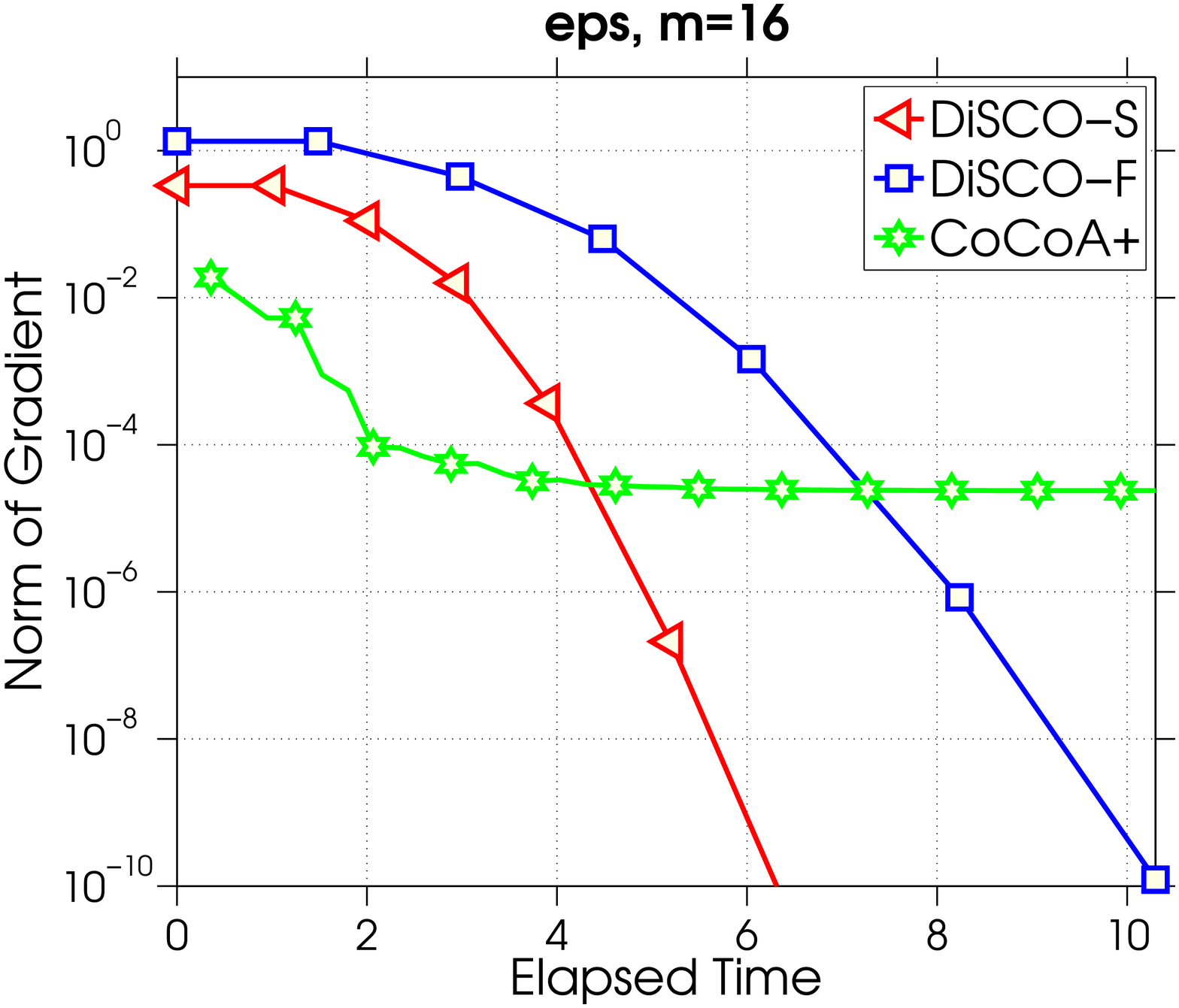}
\includegraphics[scale=.18]{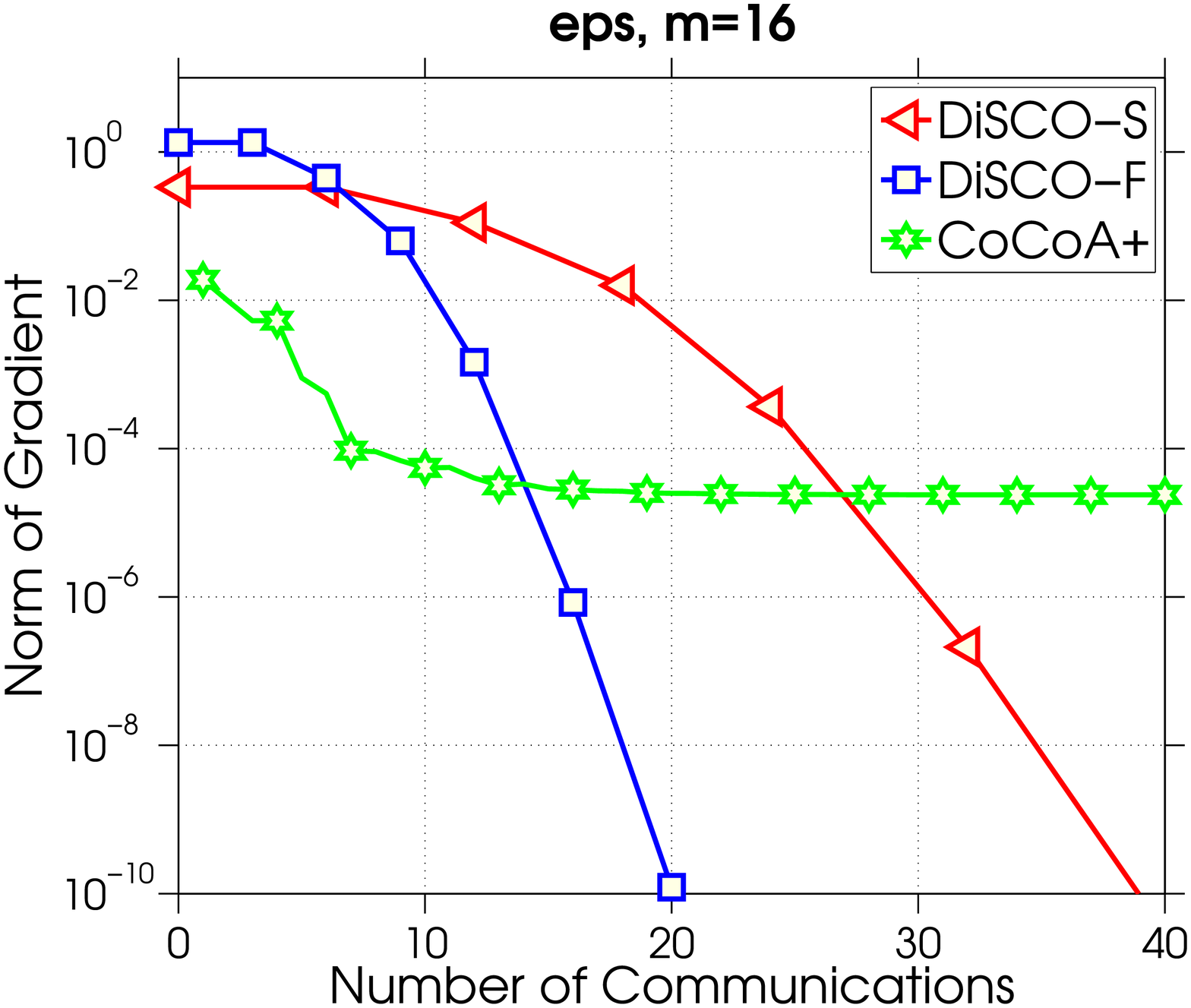}
\includegraphics[scale=.18]{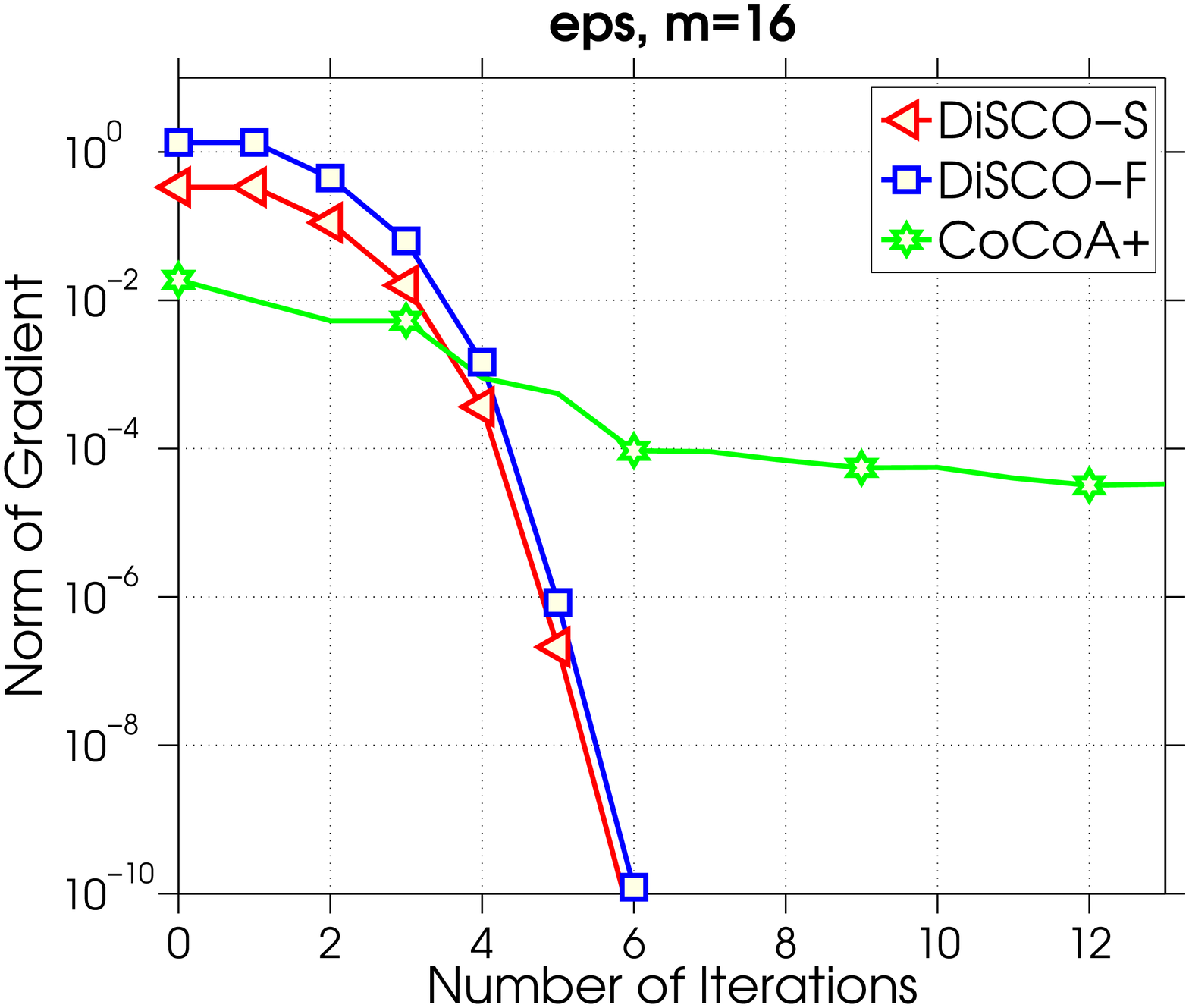}

\caption{Comparison of DiSCO-S, DiSCO-F and CoCoA+ on various datasets.
} 
\label{fig:compareRcv}
\end{figure}

We present experiments on several standard large real-world datasets: 
news20.binary $(d=1,355,191; n=19,996;  0.13GB)$;
kdd2010(test) $(d=29,890,095; n=748,401; 0.19GB)$; and 
    epsilon $(d=2,000; n=100,000;  3.04GB)$. Each data was split into $m$ machines. 
 We implement DiSCO-S, DiSCO-F and CoCoA+ \cite{ma2015adding} algorithms for comparison in C++, and run them on the Amazon cloud, using 4 m3.xlarge EC2 instances. 
Figure \ref{fig:compareRcv} compares the evolution of $\|\nabla f(w)\|$
as function of elapsed time, number of communications and iterations.
As it can be observed, the DiSCO-F needs almost the same number of iterations as DiSCO-S, however, it needs roughly just half the communication, therefore it is much faster (if we care about elapsed time).

\clearpage 
 \bibliographystyle{plain} 
 
\bibliography{literature}

\end{document}